\documentclass{article}
\newcommand{\R}{I\!\! R}

\begin{document}

Tsemo Aristide

103-135 Leeward Glenway

North York

M3C 2Z6

tsemoaristide@hotmail.com

\bigskip

\bigskip

{\bf KAHLER, POISSON GEOMETRY OF CR LIE GROUPS.}

\bigskip

\bigskip
                     \centerline{\it "Las frutas de la honestidad se
                     recogen}

                    \centerline {\it en muy poco tempo y duran para
                    siempre."}

                              \centerline{\it Colombian Proverb.}

\bigskip

\bigskip

\centerline{\bf Abstract.}

A Cauchy Riemann  $(CR)$ Lie group is a Lie group $G$ which Lie
algebra ${\cal G}$ has a  vector subspace ${\cal H}$ endowed with
an endomorphism $j$ such that $j^2=-Id$, and for each elements
$x$, $y$ in ${\cal H}$, we have $[j(x),j(y)]-[x,y]$ is an element
of ${\cal H}$, and $[j(x),j(y)]-[x,y]=j([j(x),y] + [x,j(y)])$. In
this paper we study the geometry of  a CR Lie groups $G$ when its
Lie algebra ${\cal G}$ is endowed with more geometric structures
compatible with $j$, as kahler, and poisson type structures.

\bigskip

\bigskip

\centerline{\bf 0. Introduction.}

\bigskip

Let $(M,j)$ be a complex manifold, and $H$ an hypersurface of $M$,
for each element $x$ of $H$ the tangent space $TH_x$ of $H$ at $x$
is endowed with a maximal complex vector space $E_x=TH_x\cap
j(TH_x)$. The collection of vector spaces $E_x$ defines a vector
bundle $E$  such that for sections $X$, and $Y$ of $E$, we have:
$$
[j(X),j(Y)]-[X,Y]\in E
$$
since $(M,j)$ is a complex structure, we also have:

$$
[j(X),j(Y)]-[X,Y] = j([j(X),Y] + [X,j(Y)]).
$$

More generally, a Cauchy Riemann (CR) structure on a manifold $M$,
is defined by a subbundle $H$ of $TM$ endowed with an endomorphism
$j$ such that:
$$
j^2 = -Id_H
$$

For each sections $X$ and $Y$ of $H$, we have:

$$
[j(X),j(Y)] - [X,Y]\in H
$$

$$
[j(X),j(Y)] - [X,Y] = j([j(X),Y] + [X,j(Y)],
$$

The notion of CR manifolds is studied by many authors see[3]

In this paper we study Left-invariant $CR-$structures on Lie
groups compatible with geometric properties as poisson, and kahler
type properties. More precisely:

\medskip

{\bf Definition 0.1.}

A kahlerian complex real Lie group $(G,{\cal H},j,<,>)$, is a Lie
group $G$ endowed with a CR structure defined by the vector
subspace ${\cal H}$, $j$ is an endomorphism of ${\cal G}$ which
image is ${\cal H}$. We suppose that the following properties are
verified:

    $1.$ $j$ preserves ${\cal H}$, the restriction of  $j^2$ to
     ${\cal H}$ is $-Id$.

    $2.  [X,Y]-[jX,jY]\in {\cal H}$ if $X,Y \in {\cal H}$

    $3. [j(X),j(Y)]=[X,Y]+j([X,j(Y)]+ [j(X),Y]) if X,Y \in {\cal H}$

Moreover we will suppose that there exists a left-invariant
riemannian metric on $G$, defined by a scalar product on ${\cal
G}$ such that
 $\omega=<,j>$  is closed. This means that $\omega$ is
antisymmetric and
$$
\omega([x,y],z)+\omega([z,x],y)+\omega([y,z],x)=0
$$
We remark that the restriction of $\omega$ to ${\cal H}$ is not
degenerated.

\bigskip

We show results related to those structures, for example, we show
that a semi-simple kahlerian CR Lie group such that the
codimension of ${\cal H}$ is $1$ is locally isomorphic to $so(3)$
or $sl(2)$.

\bigskip

{\bf 1. The structure of kahlerian Lie groups.}

 \bigskip

Let $G$ be Lie group. A left symmetric structure on $G$, is
defined by a product on its Lie algebra $({\cal G},[,])$
$$
{\cal G}\times {\cal G}\longrightarrow {\cal G}
$$

$$
(x,y)\longrightarrow xy
$$
which verified
$$
xy - yx = [x,y]
$$
and
$$
x(yz) - (xy)z = y(xz) - (yx)z.
$$
This is equivalent to endows $G$ with a left invariant connection
which curvature and torsion forms vanish identically.

\bigskip

In this part we consider a CR kahlerian Lie group $(G,{\cal
H},j,<,>)$. The Lie algebra ${\cal G}$ of $G$, will be call a CR
kahlerian Lie algebra.

We will define the following  product on ${\cal H}$:

For $x$, $y$, $z$ $\in {\cal H}$, we set
$$
\omega(xy,z)=-\omega(y,[x,z])=\omega([x,z],y)=<[x,z],j(y)>
$$

\medskip

{\bf Proposition 1.1.}

{\it For every $x$, $y$, $u$ and $z$ in ${\cal H}$ we have
$$
\omega(xy-yx,u)=\omega([x,y],u)\leqno(1)
$$
and if the bracket $[x,y]'=xy-yx$ satisfies the jacobi identity
$$
x(yz)-(xy)z=y(xz)-(yx)z \leqno (2)
$$}

\medskip

{\bf Proof.}

The proof almost copy the one of left-invariant symplectic
structures on Lie groups.

Let $u$ be an element of ${\cal H}$, for $x$, $y$, $z$ in ${\cal
H}$, we have $\omega(xy-yx,u)=-\omega(y,[x,u])+\omega(x,[y,u])$
then $(1)$ follows from the definition of $\omega$ (the closed
property).

Now we prove the second assertion:

$$
\omega(x(yz),u)=-\omega(yz,[x,u])=\omega(z,[y,[x,u]']')
$$

We also have:
$$
\omega((xy)z,u)=-\omega(z,[xy,u])=\omega(z,[u,xy])=-\omega(uz,xy)
$$
This implies that
$$
\omega(x(yz)-(xy)z,u)-\omega(y(xz)-(yx)z,u)=
$$

$$
\omega(z,[y,[x,u]]-[x,[y,u]])+\omega(uz,xy-yx)
$$
The property $(1)$ implies that:
$$
 \omega(uz,xy-yx)=\omega(uz,[x,y]')=-\omega(z,[u,[x,y]']')
 $$
 We deduce that
$$
\omega(x(yz)-(xy)z,u)-\omega(y(xz)-(yx)z,u)=
\omega(z,[y,[x,u]']'+[x,[u,y]']'+[u,[y,x]']')=0
$$

\medskip

{\bf Corollary 1.2.}

{\it Let ${\cal G},{\cal H},j)$ be a kahlerian $cr-$algebra, then
the product $(x,y)\rightarrow xy-yx=[x,y]'$ defined on ${\cal H}$
a structure of a kahlerian Lie algebra if it satisfies the Jacobi
identity.}

\medskip

{\bf Proof.}

We deduce from the property $(2)$ that the product defined on
${\cal H}$, $(x,y)\rightarrow xy$ endows ${\cal H}$ with a
structure of Left symmetric algebra which underlying Lie algebra
is $[,]'$, the morphism $j$ is defines also a complex structure on
$H$, and the scalar product $<,>$ a kahlerian structure.

\medskip

{\bf Proposition 1.3.}

{\it Consider the vector subspace ${\cal L}$ such that for each
$x\in {\cal L}$, we have $\omega(x,{\cal G})=0$, ${\cal L}$ is a
Lie subalgebra of ${\cal G}$ and is the $<,>$orthogonal vector
space of ${\cal H}$ .}

\medskip

{\bf Proof.}

Let $x$ and $y$ two elements of ${\cal L}$, for every element $z$
of ${\cal G}$ we have:
$$
\omega([x,y],z)=\omega(y,[z,x])+\omega(x,[y,z])
$$
since $x$ and $y$ are elements of ${\cal L}$, we deduce that
$\omega([x,y],z)=0$, and that $[x,y]\in {\cal L}$.

 For $x\in {\cal L}$, and $y\in {\cal
H}$, we have:
$$
\omega(x,y)=<x,j(y)>=0,
$$
Since the restriction of $j$ to ${\cal H}$ is an automorphism, we
deduce that ${\cal L}$ and ${\cal H}$ are $<,>$ orthogonal each
other

\medskip

 {\bf Proposition 1.4.}

{\it Let $L$ be the Lie subgroup which Lie algebra is ${\cal L}$,
and  $M$  the right quotient $G/L$, then $M$ is a kahlerian
manifold.}

\medskip

{\bf Proof.}

The fact that ${\cal H}$ and ${\cal L}$ are orthogonal each other
implies that the riemannian metric $<,>$ gives rise to a metric
$<,>'$ of $M$, the morphism $j$ also gives rise to a complex
structure $j'$  of $M$. Denote by $p$ the projection
$p:G\rightarrow G/L$, we have $\omega=p^*<,j'>'$. This implies
that $\omega'=<,j'>'$ is a symplectic form defined on $M$, thus
$(M,\omega',j')$ is a kahlerian manifold.

\medskip

\medskip

Let $({\cal H},[,]',j',<,>')$ be a kahlerian Lie algebra, and $V$
a vector space. Supposed defined a Lie algebra structure on ${\cal
G}={\cal H}+V$ such that there exists a map:
$$
\alpha: {\cal H}\times {\cal H}\rightarrow V
$$
such that for $x,y\in {\cal H}$ we have
$[x,y]=[x,y]'+\alpha(x,y)$, suppose that
$\alpha(j(x),j(y))=\alpha(x,y)$.

We suppose also that there exists a scalar product $<,>$ on ${\cal
G}$ which extends $<,>'$ such that ${\cal H}$ and $V$ are
orthogonal, we also extend $j'$ to an endomorphism $j$ of ${\cal
G}$ such that $j(V) =0$. We suppose that the form $\omega=<,j>$ is
closed, then
 $({\cal G},j,<,>)$ is a
$CR-$kahlerian algebra.

Remark that the fact that ${\cal G}$ is a kahlerian $CR-$Lie
algebra implies the following property:
$$
\oint \alpha([x,y]',z)+[\alpha(x,y),z]=0
$$

\medskip

Remark that if ${\cal H}$ is a sub Lie algebra of ${\cal G}$, then
its symplectic structure defined on ${\cal G}$ a left invariant
Poisson structure.

\medskip

{\bf Examples of kahlerian $CR-$structures.}

\medskip

$1.$

 Consider the  $n-$dimensional commutative Lie algebra
${\R}^n$ endowed with its flat riemannian metric $<,>$, and $V$ an
even dimension subspace of ${\R}^n$ endowed with a linear map $j$
such that $j^2=-id$, then $({\R}^n,V,j,<,>)$ is a Lie kahlerian
$CR-$algebra.

$2.$

Consider the semi-simple algebra $so(3)$, and $(e_1,e_2,e_3)$ its
basis in which its Lie structure is defined by
$$
[e_1,e_2]=e_3, [e_1,e_3]=-e_2, [e_2,e_3]=e_1
$$
On $Vect(e_1,e_2)$ the subspace generated by $e_1$ and $e_2$, we
consider the linear map $j$ defines by $j(e_1)=e_2$ and
$j(e_2)=-e_1$.

Let $<,>$ be the scalar product defined on $so(3)$ by
$<e_i,e_j>=\delta_{ij}$,

The family $(so(3),V,j,<,>)$ is a $cr-$kahlerian Lie algebra.

\medskip

{\bf Proposition 1.6.}

{\it Let $Z({\cal G})$ be the center of ${\cal G}$, then
$U=(Z({\cal G})\cap {\cal H})+ j(Z({\cal G})\cap {\cal H})$ is a
Lie commutative algebra and for every element $z$ of $U$,
$adz({\cal H})\subset {\cal H}$.}

{\bf Proof.}

Let $z$ and $z'$ be elements of $Z({\cal G})\cap {\cal H}$, we
have $[j(z),j(z')]=[z,z']+j([z,j(z')]+[j(z),z'])$. Since $z$ and
$z'$ are in the center of ${\cal G}$, we deduce that
$[j(z),j(z')]=0$.

The fact that $[z,{\cal H}]\subset {\cal H}$ follows from the fact
that for $x$, $y$ in ${\cal H}$, we have $[x,j(y)]+[j(x),y]$ is an
element of ${\cal H}$.

\medskip

\medskip

{\bf Proposition 1.7.}

{\it Let ${\cal G}$ be a $CR-$algebra, suppose that there exists
an ideal $I$ supplementary to ${\cal H}$, then ${\cal H}$ is
endowed with a complex Lie structure.}

\medskip

{\bf Proof.}

The projection $p:{\cal G}\rightarrow {\cal H}$ parallel to $I$
defines on ${\cal H}$ a Lie-complex structure.

\medskip

Conversely suppose given an extension
$$
0\longrightarrow I\longrightarrow {\cal G}\longrightarrow {\cal
U}\longrightarrow 0
$$
where ${\cal U}$ is a Lie algebra endowed with a complex
structure, then a supplementary space ${\cal H}$ of $I$ defines a
$CR-$structure on ${\cal G}$ if and only if for every $x,y$ in
${\cal H}$, $[jx,jy]-[x,y]$ is an element of ${\cal H}$ and
$[j(x),j(y)]-[x,y]=j([j(x),y]+[x,j(y)])$, where $j$ is the
pulls-back of the complex structure of ${\cal U}$.

\medskip
{\bf Proposition 1.8.}

{\it Let $(G,{\cal H}<,>,j)$ be a CR kahlerian Lie group. Suppose
that ${\cal H}^{\perp}= ker j$ is an ideal, then ${\cal G}$ is not
semi-simple.}

\medskip

{\bf Proof.}

We have seen that if ${\cal H}^{\perp}$ then ${\cal H}$ is endowed
with a structure of a kahlerian Lie algebra which is known not be
semi-simple. Since the quotient of ${\cal G}$ by ${\cal
H}^{\perp}$ is isomorphic as a Lie algebra to ${\cal H}$, we
deduce that ${\cal G}$ is not semi-simple.

\medskip

{\bf Proposition 1.9.}

{\it Suppose that ${\cal H}^{\perp}$ is an ideal, then $G$ is the
product of a family of groups $G_i$, where $G_0$ is flat, and for
$i\geq 1$, the holonomy of the riemannian structure of $G_i$ is
irreductible. Each group $G_i$ is endowed with a kahlerian
$CR-$structure defined by the subspace ${\cal H}_i$ of the Lie
algebra ${\cal G}_i$ of $G_i$ such that the sum of the dimension
of ${\cal H}_i$ is ${\cal H}$.}

\medskip

{\bf Proof.}

Suppose that ${\cal H}^{\perp}$ is an ideal, then the quotient
${\cal L}$ of ${\cal G}$ by ${\cal H}^{\perp}$ is a kahlerian Lie
algebra.   The theorem (Lichnerowicz Medina [12]), implies that
${\cal L}=\sum {\cal L}_i$ where each ${\cal L}_i$ is a kahlerian
Lie algebra. Now consider the De Rham decomposition of $G$ as a
product of groups $G_i$. The projection $G\rightarrow L$ respect
this decomposition. Suppose that the projection $p:{\cal
G}_i\rightarrow {\cal L}_i$ is not trivial, then the orthogonal of
the kernel of $p$ defines ${\cal H}_i$.

\medskip

Consider the set of functions $C^{\infty}(G_H)$ defined on $G$
such that $d^nf\in S^n(T^*{\cal H})$ where $T{\cal H}$ is the left
invariant distribution on $G$ defined by ${\cal H}$. For each
$f\in C^{\infty}(G_H)$, there exists a vector field $X_f$ such
that
$$
\omega(X_f,.)=df
$$
we will denote by $\{f,g\}=\omega(X_f,X_g)$.

\medskip

{\bf Proposition 1.10.}

{\it The algebra $(C^{\infty}(G_H),\{,\}$ is a Poisson algebra,
i.e $\{,\}$  verifies the Jacobi identity.}

\medskip

\bigskip

{\bf 2. Homogeneous Kahler CR manifold.}

\bigskip

Let $(G,{\cal H},j,<,>)$ be a kahler CR Lie group and $\Gamma$ a
cocompact discrete subgroup of $G$, the manifold $M= G/\Gamma$
inherits a CR structure from  $G$.

The orbits of the left action of the group $L$ on $G$ defines a
foliation  ${\cal F}_L$ on $M$. (The Lie algebra of $L$ is ${\cal
H}$ $<,>$ orthogonal.

\medskip

{\bf Proposition 2.1.}

{\it The orbits of the foliation ${\cal F}_L$ are closed if and
only if   the group $\Gamma L$ is closed in $G$, in this case $M$
is the total space of a fibration over a kahler manifold.}

\medskip

{\bf Proof.}

Suppose that the group $\Gamma L$ is closed, then the quotient $G/
\Gamma L$ is a Kahlerian manifold $N$, the projection map
$M\rightarrow N$ induced by the identity map of $G$ which is the
given fibration.

Conversely suppose that the orbits of the foliation ${\cal F}_L$
are closed. Consider a sequence $g_n$ of $\Gamma L$ which
converges towards the element $g$ of $G$. Consider a neighbourhood
$U$ of $g$ such that the restriction of the projection $p:
G\rightarrow M$ to $U$ is injective, then $p(g)$ is an element of
the adherence of $p((\Gamma L)e)={\cal F}_{p(e)}$, where $e$ is
the neutral element of $G$. Since we have supposed that the leaves
of ${\cal F}_L$ are closed, $p(g)$ is an element of ${\cal
F}_{p(e)}$ which means that $\Gamma L$ is closed in $G$. We can
thus define the Kahler manifold $G/\Gamma L$.

\medskip

{\bf Deformation of CR kahlerian structures of homogeneous
manifolds.}

\medskip

Let $(G,{\cal H},j,<,>)$ be a Lie group endowed with a CR
kahlerian structure. Consider a cocompact subgroup $\Gamma$ of $G$
the manifold $M=G/\Gamma$ inherits from $G$ a CR structure. In
this section, we will define the deformation of those structures
from two points of view.

\medskip

Supposed fixed the CR kahlerian structure of $G$, and a compact
manifold $M$, Let $\Gamma$ be  a group we consider the set of
representations $R(\Gamma,G)$ such that for each $u\in
R(\Gamma,G)$, $u$ is injective and $G/u(\Gamma)$ is a compact
manifold.

To elements $u$ and $u'$ of $R(\Gamma,G)$ will be said equivalent
if and only there exists an element $g$ of $G$ such that
$u'=gug^{-1}$. We denote by $Def_1(\Gamma,G,{\cal H},j,<,>)$ the
space of equivalence classes of those CR kahlerian structures.

\medskip

Now consider $RCK(G)$, the set of real complex kahlerian
structures of $G$, then for a cocompact subgroup $\Gamma$ of $G$
$M=G/\Gamma$ inherits a CR kahlerian structure for each element
$u$ of $RCK(G)$ denotes by $(M,u)$. We will say that $(M,u_1)$ is
equivalent to $(M,u_2)$ if and only if there exists an isomorphism
of CR kahlerian complex manifolds between $(M,u_1)$ and $(M,u_2)$.
and denote by $Def_2(M, G)$ the set of those CR structures.

\bigskip

{\bf 3. Kahlerian  codimension $1$ $CR-$structures.}

\medskip

{\bf Theorem 3.1.}

{\it Suppose that the codimension of ${\cal H}$ is $l$, and $G$ is
semi-simple then
 $G$ is a Lie group of rank $\leq l$.}

 \medskip

 {\bf Proof.}

Let $({\cal G}, {\cal H},j,<,>)$ be a Lie semi-simple algebra
endowed with a codimension $l$ $CR-$kahlerian structure. This
means that the codimension of ${\cal H}$ in ${\cal G}$ is $l$.

The map ${\cal G}\rightarrow {\cal G}^*$,
$$
X\longrightarrow \omega(X,.)
$$
is a $1-$cocycle for the coadjoint representation. Since ${\cal
G}$ is semi-simple, this cocycle is trivial. There  exists an
element $\alpha$ of ${\cal G}^*$ such that for $x,y\in {\cal G}$
we have:
$$
\omega(x,y)=\alpha([x,y])
$$
Let $K$ be the Killing form of ${\cal G}$, there exists $X\in
{\cal G}$ such that for each $Y\in {\cal G}$ we have:
$$
K(X,Y)=\alpha(Y)
$$
The Lie algebra ${\cal L}=\{x,:\omega(x,{\cal G})=0\}$ is a
dimension $l$ subalgebra of ${\cal G}$ since the codimension of
${\cal H}$ is $l$. The Lie algebra ${\cal L}$ is the Lie algebra
of the subgroup $L$ of $G$ which preserves $\alpha$, $L$ is also
the subgroup which preserves $X$ since $K$ is invariant by the
adjoint representation. We deduce that the  rank of $G$ is less or
equal than $l$, and then that $G$ is isomorphic to $sl(2)$ or
$so(3)$ if the codimension of ${\cal H}$ is $1$.

\bigskip

{\bf Corollary 3.2.}

{\it Suppose that ${\cal G},{\cal H},<,>,j)$ is a semi-simple
codimension $1$ $cr-$structure (the codimension of ${\cal H}$ is
$1$), then
 ${\cal G}$ is $so(3)$ or $sl(2)$.}

\bigskip

{\bf 4. Poisson $CR-$structures.}

\medskip

Let $M$ be a manifold, and $TH$ and $TU$ two supplementary
subbundles of its tangent bundle $TM$.

\medskip

{\bf Definition 4.1.}

An $(TH,TU)-$pseudo-Poisson structure on $M$ is defined by a
bivector $\Lambda\in \Lambda^2TM$ such that
$$
[\Lambda,\Lambda]\in TU\Lambda^2 TM
$$

where $[\Lambda,\Lambda]$ is the schouten product of $\Lambda$ by
$\Lambda$. The bivector $\Lambda$ defined on $C^{\infty}(M)$ the
bracket $\{,\}$ by the formula
$$
\{f,g\}=\Lambda(df,dg)
$$

A morphism $f:(M,\Lambda)\rightarrow(M',\Lambda')$ is  a
differentiable map which commutes with $\{,\}$.

Suppose that the distribution $TH$ defines on $M$ a
$cr-$structure, we will say that $(M,TH,TU,j)$ defines a
pseudo-Poisson $cr-$structure on $M$, if $j$ preserves $\Lambda$.

\medskip

 Let $G$ be a Lie group which Lie algebra is ${\cal G}$.
Consider a subspace ${H}$ of ${\cal G}$ and $U$ a supplementary
space to $H$. The vector spaces $H$ and $U$ define   on $G$  right
invariant distributions $TH$ and $TU$.

\medskip

{\bf Definition 4.2.}

A pseudo-Lie Poisson structure on $G$ is a bivector $\Lambda\in
\Lambda^2TG$ such that
$$
[\Lambda,\Lambda]\in TU\Lambda^2 TG.
$$
Moreover we suppose that $\Lambda $ is multiplicative i.e that the
product $G\times G\rightarrow G$ is a morphism of pseudo-Poisson
structures.

\medskip

The bracket $\{,\}$ defined by $\Lambda$ satisfies the following
properties:
$$
\{f,f'\}(xy)=\{f\circ L_x,f'\circ L_x\}(y)+\{f\circ R_y,f'\circ
R_y\}(x).
$$

If we denote by $T_uL_x$ the differential of $L_x$ in $u$, we have

$$
\Lambda(xy)=T_yL_x\Lambda_y + T_xR_y\Lambda_x
$$

\medskip

Consider the tensors $\Lambda_R(x)=T_xR_{x_{-1}}\Lambda_x$ and
$\Lambda_L(x)=T_xL_{x_{-1}}\Lambda_x$

\medskip

{\bf Proposition 4.3.}

{\it The fact that $\pi $ is multiplicative is equivalent to

1. $\pi_R$ is a $1-$cocycle for the adjoint representation
$G\rightarrow \Lambda^2{\cal G}$,i.e $\pi_R(xy)=\pi_R(x) +
Ad_x(\pi_R(y))$.

2. $\pi_L$ is a $1-$cocycle for the adjoint action of the opposite
group of $G$ in $\Lambda^2G$ i.e $\pi_L(xy)=\pi_L(y)+
Ad_{y^{-1}}(\pi_L(x))$.}

\bigskip

Let $r$ be an element of $\Lambda^2{\cal G}$, $r_{-}$, and $r_{+}$
the left and right invariant tensors defined by $r$. We denote by
$\pi$ the tensor  $r_{+}-r_{-}$.

The tensor $\pi$ defines a pseudo-Poisson structure if and only if
$[\pi,\pi]=[r,r]_{+}-[r,r]_{-}\in TU\Lambda^2TG$, or equivalently
if
$$
Ad_x[r,r]-[r,r]\in U\Lambda^2{\cal G}.
$$

 Suppose that $H$ defines on $G$ a $cr-$structure.
We will say that $(G,H,\Lambda,j)$ is  a $cr-$Poisson structure if
$j$ preserves $\Lambda$.

\medskip

 Moreover we assume that $\Lambda$ is invariant by $j$.

\medskip

{\bf Remark.}

Suppose defined the $cr-$Poisson structure $({\cal G}_1,{\cal
H}_1,j_1,\Lambda_1)$ and $({\cal G}_2,{\cal H}_2,j_2,\Lambda_2)$,
then the tensor $\Lambda_1\times\Lambda_2$ defines a
$cr-$structure $({\cal G}_1\times{\cal G}_2,{\cal H}_1\times{\cal
H}_2,\Lambda\times\Lambda_2,j_1\times j_2$ called the product of
the Poisson $cr-$structures.

\bigskip

\centerline{\bf References.}

\bigskip

1. Aubert, A. Structures affines et pseudo-metriques invariantes a
gauche sur des groupes de Lie.

Thesis Universite de Montpellier II 1996

\medskip

2. Benson, C. Gordon C. Kahler and symplectic structures on
nilmanifolds.

Topology 27 (1988) 513-519

\medskip

3. Bland, J. Epstein, C. Embeddable CR structures and deformation
of pseudo-convex surfaces.

J. Algebraic Geometry 5 (1996) 277-368

\medskip

4. Boothby, W. Wang, H. On contact manifolds.

 Ann. Math. 68. (1958)
721-734

\medskip

5. Dardie J-M, Groupes de Lie symplectiques ou kahleriens et
double extension.

Thesis Universite de Montpellier II 1993

\medskip

6. Diatta, A. Geometrie de Poisson et de Contact des espaces
homogenes

Thesis Universite de Montpellier II 2000

\medskip

7. Domini, S. Gigante, G. Classification of left invariant CR
structures on $Gl(3,{\R})$.

 Rivist. Mat. Univ di Parma.

\medskip

8. Dorfmeister, J. Nakajima, K. The fundamental conjecture for
homogeneous kahler manifolds.

Acta Mathmatica 161 (1986) 189-208

\medskip

9. Gauthier, D. Lie CR algebras.

 Diff. Geometry and its applications
(1997) 365-376

\medskip

10. Gigante, G. Tomassini. CR structures on a real Lie algebra.

Adv. Math (1992) 67-81

\medskip

 11. Lichnerowicz, A. Les groupes Kahleriens,

  Symplectic
geometry and mathematical physics (Edit P. Donato et al). Pro.
Math 99

\medskip

12. Lichnerowicz A, Medina, A. On Lie groups with left invariant
symplectic or kahlerian structures.

Lett. Math. Phs. 16 (1988)225-235

\medskip

13. Milnor, J. Curvature of left invariant metric on Lie groups.

Adv. Math. (1976) 293-329

\medskip

14. Tsemo, A. Kahlerian, Poisson CR structures on Lie groups and
the method of double extension. In preparation.

\end{document}